\documentclass[12pt,a4paper,reqno]{amsart}

\usepackage{anysize}
\usepackage[utf8]{inputenc}
\usepackage{amsmath}
\usepackage{amssymb}
\usepackage{array}

\marginsize{2.5cm}{2.5cm}{3cm}{3cm}

\newtheorem{example}{Example}

\begin{document}

\title{On $p$-adic expansions of Ramanujan-like series}
\author{Jesús Guillera}
\address{Department of Mathematics, University of Zaragoza, 50009 Zaragoza, SPAIN}
\email{jguillera@gmail.com}

\maketitle

\begin{abstract}
Inspired by a Zudilin-Zhao's supercongruences pattern related to Ramanujan-like series for $1/\pi^k$, we conjecture a kind of $p$-adic expansions.
\end{abstract}

\section{Introduction}

Zudilin \cite{Zu0, GuiZu} observed a pattern of supercongruences related to Ramanujan series for $1/\pi^k$, and Zhao \cite{Zhao-2} added a term to it. Here is an example
\begin{equation}\label{p-adic-exam-01}
\sum_{n=0}^{p-1}  \frac{\left( \frac12 \right)_n^5}{(1)_n^5} \frac{(-1)^n}{4^n} (20n^2+8n+1) \, \, {\overset{?} \equiv} \, \, p^2 - \frac72 \zeta(4-p) \, p^5 \pmod{p^6}
\end{equation}
for primes $p>3$. They are supercongruences related to the Ramanujan-like series
\begin{equation}\label{gui-20}
\sum_{n=0}^{\infty}  \frac{\left( \frac12 \right)_n^5}{(1)_n^5} \frac{(-1)^n}{4^n} (20n^2+8n+1) = \frac{8}{\pi^2},
\end{equation}
proved by the author \cite{GuiEMpi2}. In this paper we try to find out where results like (\ref{p-adic-exam-01}) could come from.

\section{An heuristic argument} 

If we expand (\ref{gui-20}) extended with $x$ in powers of $x$, we get 
\begin{multline}\label{rama-x}
\sum_{n=0}^{\infty}  \frac{\left( \frac12 \right)_{n+x}^5} {(1)_{n+x}^5} (-1)^n \frac{20(n+x)^2+8(n+x)+1}{4^{n+x}} \\ 
{\overset{?} =} \, \, \frac{8}{\pi^2} - 4 x^2 + 50 \zeta(2) x^4 - 448 \zeta(3) x^5 + O(x^6).
\end{multline}
A $p$-adic analogue of $1/\pi^2$ is $p^2$, and a $p$-adic analogue of $\zeta(k)$ is $\zeta_p(k)$, which is defined in \cite[Definition 3.4]{Rosen}. Then, we conjecture that a $p$-adic analogue of the identity (\ref{rama-x}) is an expansion (with rational coefficients) of the form:
\begin{equation}\label{p-adic-01}
\sum_{n=0}^{p-1}  \frac{\left( \frac12 \right)_n^5}{(1)_n^5} \frac{(-1)^n}{4^n} (20n^2+8n+1) 
\, \, {\overset{?} \equiv} \, \, r_2 p^2 + r_4 \zeta_p(2) \, p^4  + r_5 \zeta_p(3) \, p^5  + \cdots.
\end{equation}
But $\zeta_p(2)=0$ \cite[Section 1.3]{Rosen-2} and \cite[answer by Loeffler]{ustinov}, and $\zeta_p(3) \equiv \zeta(4-p) \pmod{p}$ for $p \geq 5$ is a particular case of  $\zeta_p(k) \equiv \zeta(1+k-p) \pmod{p}$ for $p \geq k+2$ (which comes from \cite[page 10 top]{Rosen}). Hence, we see that (\ref{p-adic-01}) agrees with (\ref{p-adic-exam-01}), and we conjecture the $p$-adic expansion identity:
\begin{equation}\label{p-adic-02}
\sum_{n=0}^{p-1}  \frac{\left( \frac12 \right)_n^5}{(1)_n^5} \frac{(-1)^n}{4^n} (20n^2+8n+1) 
\, \, {\overset{?} \equiv} \, \, p^2 - \frac72 \zeta_p(3) \, p^5  + \cdots.
\end{equation}

\section{Other examples}

\begin{example} \rm
The following formula holds \cite[Example 60]{Chu-Zhang}:
\begin{equation}\label{otro}
\sum_{n=0}^{\infty}  \frac{\left( \frac12 \right)_n \left( \frac16 \right)_n \left(\frac56 \right)_n}{(1)_n^3} \frac{74n+7}{2n+1} \left(\frac{27}{64}\right)^n = 8.
\end{equation}
Expanding in powers of $x$, we get
\begin{multline}\label{otro-x}
\sum_{n=0}^{\infty}  \frac{\left( \frac12 \right)_{n+x} \left( \frac16 \right)_{n+x} \left(\frac56 \right)_{n+x}}{(1)_{n+x}^3} \frac{74(n+x)+7}{2(n+x)+7} \left(\frac{27}{64}\right)^{n+x} \\ 
\, \, {\overset{?} =} \, \, 8 -144 \zeta(2) x^2 + 1792 \zeta(3) x^3 +O(x^4).
\end{multline}
Hence, we conjecture the following $p$-adic expansion:
\[
\sum_{n=0}^{p-1}  \frac{\left( \frac12 \right)_n \left( \frac16 \right)_n \left(\frac56 \right)_n}{(1)_n^3} \frac{74n+7}{2n+1} \left(\frac{27}{64}\right)^n \, {\overset{?} \equiv} \, \, r_1 +r_2 \zeta_p(2) p^2 + r_3 \zeta_p(3) p^3 +\cdots.
\]
As $\zeta_p(2)=0$ and $\zeta_p(3) \, \equiv \, \zeta(4-p)$ for $p>3$, it agrees with Zudilin-Zhao's pattern of supercongruences:
\begin{equation}
\sum_{n=0}^{p-1}  \frac{\left( \frac12 \right)_n \left( \frac16 \right)_n \left(\frac56 \right)_n}{(1)_n^3} \frac{74n+7}{2n+1} \left(\frac{27}{64}\right)^n \, {\overset{?} \equiv} \, \, 7 - \frac{105}{2} \zeta(4-p) p^3  \pmod{p^4},
\end{equation}
and we conjecture the following $p$-adic expansion:
\begin{equation}
\sum_{n=0}^{p-1}  \frac{\left( \frac12 \right)_n \left( \frac16 \right)_n \left(\frac56 \right)_n}{(1)_n^3} \frac{74n+7}{2n+1} \left(\frac{27}{64}\right)^n \, {\overset{?} \equiv} \, \, 7 - \frac{105}{2} \zeta_p(3) p^3 + \cdots.
\end{equation}
\end{example}

\begin{example} \rm
The following conjectured formula is known \cite{GuiEMpi2}:
\begin{equation}\label{gui-5418}
\sum_{n=0}^{\infty}  \frac{\left( \frac12 \right)_n \left( \frac13 \right)_n \left( \frac23 \right)_n \left( \frac16 \right)_n \left(\frac56 \right)_n}{(1)_n^5} \frac{(-1)^n}{80^{3n}} (5418n^2+693n+29)\, {\overset{?} =} \, \frac{128 \sqrt 5}{\pi^2},
\end{equation}
If we expand (\ref{gui-5418}) extended with $x$ in powers of $x$, we get 
\begin{multline}\label{rama-5418-x}
\frac{1}{128} \sum_{n=0}^{\infty}  \frac{\left( \frac12 \right)_{n+x} \left( \frac13 \right)_{n+x} \left( \frac23 \right)_{n+x} \left( \frac16 \right)_{n+x} \left(\frac56 \right)_{n+x}}{(1)_{n+x}^5} (-1)^n \frac{5418(n+x)^2+693(n+x)+29}{80^{3(n+x)}} \\ 
\, \, {\overset{?} =} \, \, \frac{\sqrt{5}}{\pi^2} - \frac{15}{2} \sqrt{5} x^2 + \frac{110875}{32} L_5(2) x^4 - 42000 L_5(3) x^5 + O(x^6).
\end{multline}
We conjecture a $p$-adic analogue of the form:
\begin{multline}\label{p-adic-03}
\sum_{n=0}^{p-1}  \frac{\left( \frac12 \right)_n \left( \frac13 \right)_n \left( \frac23 \right)_n \left( \frac16 \right)_n \left(\frac56 \right)_n}{(1)_n^5} \frac{(-1)^n}{80^{3n}} (5418n^2+693n+29) \\
\\
{\overset{?} \equiv} \,
r_2 \left( \frac{5}{p} \right) \, p^2 + r_4 L_{5,p}(2) \, p^4  + r_5 L_{5,p}(3) \, p^5  + \cdots,
\end{multline}
where $(\frac{5}{p})$ is the $p$-adic analogue of $\sqrt 5$, $p^2$ the $p$-adic analogue of $1/\pi^2$, and $L_{5,p}(k)$ the $p$-adic anologue of $L_5(k)$. As  $L_{5,p}(2)=0$ and $L_{5,p}(k) \equiv L_{5}(1+k-p) \pmod{p}$ for $p \geq k+2$, it agrees with  Zudilin-Zhao's pattern of supercongruences. We get
\begin{multline}\label{p-adic-exam-04}
\sum_{n=0}^{p-1}  \frac{\left( \frac12 \right)_n \left( \frac13 \right)_n \left( \frac23 \right)_n \left( \frac16 \right)_n \left(\frac56 \right)_n}{(1)_n^5} \frac{(-1)^n}{80^{3n}} (5418n^2+693n+29) \\ {\overset{?} \equiv} \, \,  29 \left(\frac{5}{p}\right) \, p^2 - \frac{35}{216} L_{5}(4-p) \, p^3 \pmod{p^4},
\end{multline}
and we conjecture the following $p$-adic expansion
\[
\sum_{n=0}^{p-1}  \frac{\left( \frac12 \right)_n \left( \frac13 \right)_n \left( \frac23 \right)_n \left( \frac16 \right)_n \left(\frac56 \right)_n}{(1)_n^5} \frac{(-1)^n}{80^{3n}} (5418n^2+693n+29) \, \, {\overset{?} \equiv} \, \, 29 \left(\frac{5}{p}\right) \, p^2 - \frac{35}{216} L_{5,p}(3) \, p^3 + \cdots.
\]
\end{example}

\begin{example} \rm
The following Ramnujan-like series was discovered by Boris Gourevitch \cite{GuiEMpi2}:
\begin{equation}\label{gou-168}
\sum_{n=0}^{\infty}  \frac{\left( \frac12 \right)_n^7}{(1)_n^7} \left( \frac{1}{64}\right)^n (168 n^3+76 n^2+14n+1)\, {\overset{?} =} \, \frac{32}{\pi^3}.
\end{equation}
Extending it with a variable $x$ and expanding in powers of $x$, we get:
\begin{multline}\label{gou-168-x}
\frac{1}{32} \sum_{n=0}^{\infty}  \frac{\left( \frac12 \right)_{n+x}^7}{(1)_{n+x}^7} \left( \frac{1}{64}\right)^{n+x} (168 (n+x)^3+76 (n+x)^2+14(n+x)+1) \\ 
{\overset{?} =} \, \frac{1}{\pi^3} - \frac{1}{\pi} \, x^2 + \frac{16}{3} L_{-4}(1) \, x^4 - \frac{8224}{45}L_{-4}(3) \, x^6 + 1536 L_{-4}(4) \, x^7 + O(x^8).
\end{multline}
We conjecture a $p$-adic analogue of the following form:
\begin{multline}\label{gou-p-adic}
\sum_{n=0}^{p-1}  \frac{\left( \frac12 \right)_n^7}{(1)_n^7} \left( \frac{1}{64}\right)^n (168 n^3+76 n^2+14 n+1) \\ 
{\overset{?} \equiv} \, r_3 \left( \frac{-4}{p} \right) p^3 + r_4 L_{-4, p}(1) \, p^4 + r_6 L_{-4, p}(3) \, p^6 + r_7 L_{-4, p}(4) \, p^7  \pmod{p^8}.
\end{multline}
As  $L_{-4,p}(1)=L_{-4,p}(3)=0$, and $L_{-4,p}(4) \equiv L_{-4}(5-p) \pmod{p}$, we recover Zhao's type congruences:
\begin{equation}\label{gou-zhao-cong}
\sum_{n=0}^{p-1}  \frac{\left( \frac12 \right)_n^7}{(1)_n^7} \left( \frac{1}{64}\right)^n (168 n ^3+76 n^2+14 n+1) \, {\overset{?} \equiv} \, \left( \frac{-4}{p} \right) p^3 - 6 L_{-4}(5-p) \, p^7 \pmod{p^8},
\end{equation}
and we have the $p$-adic expansion
\[
\sum_{n=0}^{p-1}  \frac{\left( \frac12 \right)_n^7}{(1)_n^7} \left( \frac{1}{64}\right)^n (168 n ^3+76 n^2+14 n+1) \, {\overset{?} \equiv} \, \left( \frac{-4}{p} \right) p^3 -6 L_{-4,p}(4) + \cdots.
\]
\end{example}
\begin{example} \rm
For the formula of Ramanujan-Orr type 
\begin{equation}\label{fam-orr-ex-17}
\frac{3}{529} \sum_{n=0}^{\infty}  \frac{\left( \frac18 \right)_n \left( \frac38 \right)_n \left( \frac58 \right)_n \left(\frac78 \right)_n}{\left( \frac12 \right)_n (1)_n^3} \frac{6970n^2+4037n+280}{2n+1} (-1)^n \left(\frac{2^{14}}{23^4}\right)^n = \frac{\sqrt{23}}{\pi},
\end{equation}
proved in \cite{Gui-family-orr}, we have
\begin{multline}\label{fam-orr-ex-17-x}
\frac{3}{529} \sum_{n=0}^{\infty}  \frac{\left( \frac18 \right)_{n+x} \left( \frac38 \right)_{n+x} \left( \frac58 \right)_{n+x} \left(\frac78 \right)_{n+x}}{\left( \frac12 \right)_{n+x} (1)_{n+x}^3} \frac{6970(n+x)^2+4037(n+x)+280}{2(n+x)+1} (-1)^n \left(\frac{2^{14}}{23^4}\right)^{n+x}  \\ {\overset{?}=} \,  \frac{\sqrt{23}}{\pi} - \frac{69}{2} L_{-23}(1) x^2 + 529 L_{-23}(2)x^3 + O(x^4).
\end{multline}
The $p$-adic analogue is of the following form:
\begin{multline}\label{fam-orr-ex-17-cong}
\sum_{n=0}^{p-1}  \frac{\left( \frac18 \right)_n \left( \frac38 \right)_n \left( \frac58 \right)_n \left(\frac78 \right)_n}{\left( \frac12 \right)_n (1)_n^3} \frac{6970 n^2+4037n+280}{2n+1} (-1)^n \left(\frac{2^{14}}{23^4}\right)^n \, {\overset{?} \equiv} \\ r_1 \left( \frac{-23}{p} \right) p + r_2 L_{-23, p}(1) p^2 + r_3 L_{-23, p}(2)p^3 \pmod{p^4}.
\end{multline}But $L_{-23, p}(1)=0$ and $L_{-23, p}(2) \equiv L_{-23}(3-p) \pmod{p}$, and we see that the congruences that we have found agree with Zhao-Zudilin's pattern. We get
\begin{multline}
\sum_{n=0}^{p-1}  \frac{\left( \frac18 \right)_n \left( \frac38 \right)_n \left( \frac58 \right)_n \left(\frac78 \right)_n}{\left( \frac12 \right)_n (1)_n^3} \frac{6970 n^2+4037n+280}{2n+1} (-1)^n \left(\frac{2^{14}}{23^4}\right)^n \, {\overset{?} \equiv} \\ 280 \left( \frac{-23}{p} \right) p + 280 L_{-23, p}(2)p^3 + \cdots,
\end{multline}
and we conjecture that it is the $p$-adic analogue of (\ref{fam-orr-ex-17-x}).
\end{example}

\end{document}